\documentclass[12pt]{amsart}
\usepackage{amscd, amssymb, amsfonts, latexsym}
\usepackage{epsfig}
\usepackage{graphicx}
\newtheorem{Thm}{Theorem}
\newtheorem{Prop}{Proposition}

\newtheorem{Def/Thm}{Definition/Theorem}
\newtheorem{Cor}{Corollary}

\newtheorem{Lemma}{Lemma}

\newcommand{\ti }{\times}

\newcommand{\ra }{\rightarrow}

\newcommand{\Spec}{{\mathrm{Spec}}}

\newcommand{\PP }{{\mathbb P}}

\newcommand{\CC }{{\mathbb C}}
\newcommand{\ZZ }{{\mathbb Z}}
\newcommand{\RR }{{\mathbb R}}
\newcommand{\ke }{{\varepsilon }}
\newcommand{\kb }{{\beta}}
\newcommand{\ka }{{\alpha}}
\newcommand{\kg }{{\gamma}}

\newcommand{\kl }{{\lambda}}

\begin{document}

\author{Dosang Joe and Bumsig Kim}
\title{Equivariant mirrors and the Virasoro conjecture for
flag manifolds}
\date{Revised on October 24, 2002. To appear in International Mathematics Research Notices}
\begin{abstract}
We found an explicit description of all $GL(n,\RR)$-Whittaker functions as oscillatory integrals
and thus constructed equivariant mirrors of flag manifolds.
As a consequence we proved
the Virasoro conjecture for flag manifolds.
\end{abstract}

\maketitle

\pagestyle{plain}

\section{Introduction}

A general quintic hypersurface in $\CC P^4$ is a compact
Calabi-Yau threefold. Its rational Gromov-Witten invariants had
been predicted by mirror symmetry discovered in string theory
\cite{COGP}. The prediction was proven by Givental \cite{Ge},
which is to be explained as follows. Let $n_d$ be the virtual
number of degree $d$ rational curves in the quintic threefold and
let $F(q)=1+ \frac{1}{5}\sum _{d=1}^{\infty}n_d
d^3\frac{q^d}{1-q^d}$. The Givental $J$-function $J$ for the
quintic hypersurface, satisfies the so-called quantum
differential equation
\[ (\hbar\frac{d}{dt})^2 \frac{1}{F(q)} (\hbar\frac{d}{dt})^2 J=0\]
 where $t=\ln q$ and $\hbar$ is a
formal parameter. The prediction was that the solutions to the
quantum differential equation coincide with those to the
Picard-Fuchs differential equation of the mirror family after the
explicit mirror transformation. The periods I, the solutions to
the Picard-Fuchs differential equation, are $$I=\int _{\Gamma}
 \frac{\bigwedge _{i=1}^5 dx_i
}{d(x_1+...+x_5 -1)\wedge d(x_1....x_5-q) }$$ where
  $\Gamma$ are  real 3-cycles in the affine varieties $\{ x_1+...+x_5=1, x_1...x_5=q \}$
birational to the mirror manifolds.

In 1993 Givental proposed a generalization of this mirror
phenomenon to non-Calabi-Yau manifolds. A mirror family of a Fano
manifold $X$ is by definition a stationary phase integral
representation
$$I=\int _{\Gamma \subset Y_q} e^{f_q/\hbar} \omega _q $$
 of solutions
to the quantum differential equations of $X$ up to change of coordinates, where $q\in
H^2(X,\CC )/2\pi \sqrt{-1}H^2(X,\ZZ)$. Here $f_q$ is a
holomorphic function on a possibly noncompact variety $Y_q$ and
$\omega _q$ is a holomorphic volume form on $Y_q$. The mirror
theorem has been established for Calabi-Yau or Fano complete
intersections $X$ in a projective space \cite{Ge} and further
generalized to the case when the ambient space is a Fano toric
projective manifold \cite{Gt}.

The quantum differential operators for the flag manifolds have
been found in \cite{K, Gs} to be ``quantum Toda operators\rq\rq
which are by definition nonconstant integrals of motions for
quantum Toda lattice. Moreover a mirror family of the flag
manifold has been constructed in \cite{Gs}. In this paper we
build the equivariant mirror of the flag manifold, that is, a
stationary phase integral representation of complete spectra of
quantum Toda operators (see theorem \ref{sol}). Using the
equivariant mirror construction we confirm the so-called
``R-conjecture\rq\rq which leads to a proof of a formula for any genus
gravitational descendent potential for flag manifolds and as
a corollary we obtain the Virasoro
conjecture for the manifolds
(see theorem \ref{Rconjforflag}, theorem \ref{formula}, and
corollary \ref{virasoro}). The equivariant version of the formula was
obtained in \cite{Gsemi}.

Finally we remark that the integral representation that we found
gives another explicit description of all $GL(n,\RR)$-Whittaker
functions \cite{Ja, St}.

\section{Quantum Differential Operators for Flag Manifolds }
\subsection{Quantum Toda Operators}

Denote by $t_i$, $i=0,...,n$, the standard coordinate functions on $\CC ^{n+1}$,
and let $q_i =e^{t_i-t_{i-1}}$, $i=1,...,n$.
Consider a  matrix
\[A =\left[\begin{array}{ccccc}
p_0 & q_1 &  0  &  0  & ... \\
-1  &p_1 & q_2 &  0  & ... \\
  \  &  .  &  .  &  .  & \   \\
 0  & ... &  0  & -1  &p_n
\end{array} \right] \ \] of size $(n+1) \ti (n+1)$.
If  $$\det (A + xI) = x ^{n+1} + \sum _{i=1}^{n+1}D_i
(p_0,...,p_n, q_1,...,q_n) x ^{n+1-i},$$ denote
$D_i(\hbar\frac{\partial}{\partial t_0},...,
\hbar\frac{\partial}{\partial t_n}, q_1,...,q_n)$ simply by $D_i$
unless stated otherwise. Let \[ H = \frac{\hbar ^2}{2}\sum
_{i=0}^n\frac{\partial^2}{\partial t_i^2} - \sum
_{i=1}^ne^{t_i-t_{i-1}} \] which is a quantization of the
Hamiltonian of non-periodic Toda lattice. The quantum Toda
lattice is a completely integrable system with quantum integrals
$D_i$, $i=1,...,n+1$. We include an elementary proof of the
commutativity of $D_i$'s which is known.

\begin{Prop}
Let $D$ be a linear holomorphic
differential operator on $\CC ^{n+1}$ with coefficients in the
Laurent polynomial ring $\CC [\hbar, e^{\pm t_0},...,e^{\pm t_n}]$ over $\CC [\hbar]$.
Suppose that $[H, D]=0$, then the coefficients of the principal part of $D$
are contained in $\CC [\hbar]$.
\end{Prop}

{\em Proof.\/ }  Let $m$ be the order of the differential operator $D$ and for
multiple index $\ka$, let $\partial ^{\ka}
=\frac{\partial ^{\ka_0}}{\partial t_0  ^{\ka_0}}
... \frac{\partial ^{\ka_n}}{\partial t_n  ^{\ka_n}}$.
Then $D=\sum _{|\ka | =m} a_{\ka} \partial ^{\ka} + D'$ where $a_\ka$ are some polynomials
in $e^{\pm t_0}$,...,$e^{\pm t_n}$ over $\CC[\hbar]$
and $D'$ is a differential operator with order less than $m$.
In the bracket $[H,D]$, the degree $m+1$ part is created
only in $[\Delta,\sum _{|\ka | =m} a_{\ka} \partial ^{\ka}]$
where $\Delta$ is the Laplace operator
$\sum \frac{\partial ^2}{\partial t_i ^2}$.
It is enough to prove that if $[\Delta,\sum _{|\ka | =m} a_{\ka} \partial ^{\ka}]$
has no order $m+1$ part then $a_\ka \in \CC [\hbar]$ for all $\ka$.
To prove the claim, use induction on $n$. When $n=0$, it is clear.
Let $n>0$. Without loss of generality, we may assume that $a_\ka $
are not in $\CC[\hbar]$ for all $\ka$. Let $k$ be
the maximum of $\ka _0$ for all $\ka$.
Decompose $ \sum _{|\ka | =m} a_{\ka} \partial ^{\ka} = \sum _{\ka _0 = k}a_\ka \partial ^\ka
+ \mathrm{the\ rest}$. When $\ka _0 $ is $k$, $a_\ka$ does not have terms in $e^{\pm t_0}$
since $[\Delta,\sum _{|\ka | =m} a_{\ka} \partial ^{\ka}]$ has
no order $m+1$. If
$[\Delta , \mathrm{the\ rest} ] = \sum _{\kb _0=k} b_\kb \partial ^\kb + D''$, then
every $b_\kb$ depends on $e^{\pm t_0}$. So, $[\Delta , \sum _{\ka _0=k} a_{\ka} \partial ^{\ka}]$
has no order $m+1$ part. Now we apply the induction hypothesis to $\sum _{\ka _0=k} a_{\ka} \partial ^{\ka}$ whose
coefficients do not depend on $t_0$ and conclude that
$a_\ka \in \CC[\hbar]$ if $\ka _0=k$. The conclusion
 is contradictory to the assumption that for all $\ka$,
$a_\ka$ is not in $\CC[\hbar]$.
$\Box$

\begin{Thm}
$[D_i,D_j]=0$ for all $i=0,...,n$.
\end{Thm}
{\em Proof.\/} The commutativity of $H$ and $D_i$ is proven in
\cite{Gs}. Since $[H,[D_i,D_j]]=0$, by the above proposition, the
highest order part of $[D_i,D_j]$ has coefficients in $\CC
[\hbar]$. Now it is enough to prove the following claim. For any
multi-indices $\ka$ and $\kb$ and any $a$ and $b$ in the {\em
polynomial} ring $\CC [\hbar ,q_1,...,q_n]$, if we let $[a
\partial ^\ka , b \partial ^\kb] = \sum c_{\kg}
\partial ^\kg$, then any $c_\kg$ cannot be in $\CC[\hbar]$ unless $c_\kg =0$.
Notice that $[a\partial ^\ka, b\partial ^\kb ]  $ is
$$\sum _{\ka ' + \ka '' = \ka , \ \ka ' \ne 0}\left(\begin{array}{c}\ka \\ \ka '
\end{array}\right)a (\partial ^{\ka '} b) \partial ^{\ka '' + \kb } -
\sum _{\kb ' + \kb '' = \kb ,\  \kb ' \ne 0}\left(\begin{array}{c}\kb \\ \kb '
\end{array}\right) b (\partial ^{\kb '} a) \partial ^{\kb '' + \kb }.$$
However the polynomials $a(\partial ^{\ka '}b)$
and $b(\partial ^{\kb '}a)$ in $q_i$, $i=0,...,n$, over $\CC [\hbar]$
have no constant coefficients of $\CC[\hbar]$. So,
$c_\kg$ are zeros or in
$\CC [\hbar , q_1,...,q_n] - \CC[\hbar]$.
$\Box$

\subsection{Givental's $J$-functions for flag manifolds}

Recall that the (resp. equivariant) quantum differential operators are
defined to be operators which annihilate the (resp. equivariant)
$J$-functions. For the flag manifolds they are generated by the
quantum Toda operators.

 Let $Fl(n+1)$ be the set of all complete flags $V_1\subset
V_2 \subset ...\subset V_n\subset \CC ^{n+1}$ of subspaces of $\CC
^{n+1}$, where $\dim V_i=i$. It is the flag manifold of dimension
$n(n+1)/2$. Consider the torus $T=(\CC ^\ti)^{n+1}$ action on
$Fl(n+1)$ induced from the standard $T$ action on $\CC ^{n+1}$.
Denote by $\mathbb{V}_i$ the universal subbundle over $Fl(n+1)$
with fiber $V_i$ at point $(V_1,...,V_n)$. Let $p_i$ be the
equivariant first Chern class $c_1^T(\mathbb{V}_{i+1} /
\mathbb{V}_{i})$, $i=0,...,n$. Here $\mathbb{V}_0$ is by
definition rank 0 bundle. Denote $H^*_T(\mathrm{point},\ZZ ) = H^*(\PP
^{\infty})^{\otimes n+1}=\ZZ [ \kl _0,...,\kl_n]$. Here the
following convention is used: the $i$-th equivariant Chern class
$c^T_i(\mathbb{V}_{n+1})$ is  the $i$-th elementary symmetric
polynomial $\sigma _i$ in $\kl _0,...,\kl_n$. Then it is known
that the equivariant small quantum cohomology ring of $Fl(n+1)$
is generated by $p_i$ with relations $D_i(p,q)=\sigma _i$,
$i=1,...,n+1$.
 From now on take a torus $T$ as the subgroup $\{ (a_0,...,a_n) \ | \ \prod _i a_i =1\}$
 of $(\CC^\ti) ^{n+1}$     so that $\sigma _1 =0$.

\bigskip

Let $\overline{M}_{g,m}(Fl(n+1),d )$ be the moduli of
degree $d$ stable maps
$(f,C, x_1,...,x_m)$ to $Fl(n+1)$ from $m$ pointed
prestable genus $g$ curves $C$ and let
$ev_i$ be the evaluation map at the $i$-th marked point. Thus,
$f:C \ra Fl(n+1)$ is a morphism such that
$f_*[C]=d\in H_2(Fl(n+1),\ZZ )$ and $ev ([(f,C,x_1,...,x_m)])$ is
by definition  $f(x_i)$.
Let $\psi _i$ be the first Chern class of the universal cotangent
orbi-line bundle, whose fiber at $[(f,C, x_1,...,x_m)]$ is $T^*_{x_i}C$.
Fix a basis $\{ \phi _\ka \}$ of the free
$H^*_T(\mathrm{point})$-module $H^*_T(Fl(n+1),\ZZ )$.

Define equivariant Givental's $J$-functions
$J_{\ka}(t_0,...,t_n;\kl _0,...,\kl_n)$ to be
\[  \sum_{d\in H_2(Fl(n+1),\ZZ )} q^{d}
 \int _{\overline{M}_{0,1}(Fl(n+1),d )}
\frac{ev_1^*(\phi _\ka\wedge\exp (\sum p_i t_i /\hbar ))}{\hbar(\hbar-\psi _1)},
\]
where $q^d = \Pi _{i=1}^n q_i ^{<d , c_1(\mathbb{V}_i^*)>}$ with
$q_i = e^{t_i - t_{i-1}}$.

\begin{Thm}  (\cite{K, Gs})
\[ D_i J _\ka= \sigma _i J_\ka \] for all
$i$ and $\ka$.
\end{Thm}

\section{Equivariant mirrors
of flag manifolds}

In this section we find a stationary phase integral
representation of all solutions $I$ to $D_i I = \sigma _i I$,
where $i=1,...,n+1$. To do so, we shall make use of Givental's
construction $I'= \int _{\Gamma \subset Y_q} e^{f'_q/\hbar}
\omega $ of the mirror of the flag manifold \cite{Gs}. It
satisfies differential equations $D_i I' = 0$.
 From such $I'$, in order to build spectrum solutions,
we shall add appropriate weight factors to the phase function $f'$
and show that if $I=\int _{\Gamma \subset Y_q } e^{(f'_q + \text{ the
weight factors})/\hbar } \omega _q $, then $D_i I =\sigma _i I$,
$i=1,...,n+1$.

Introduce a graph and coordinates for vertices in the graph:

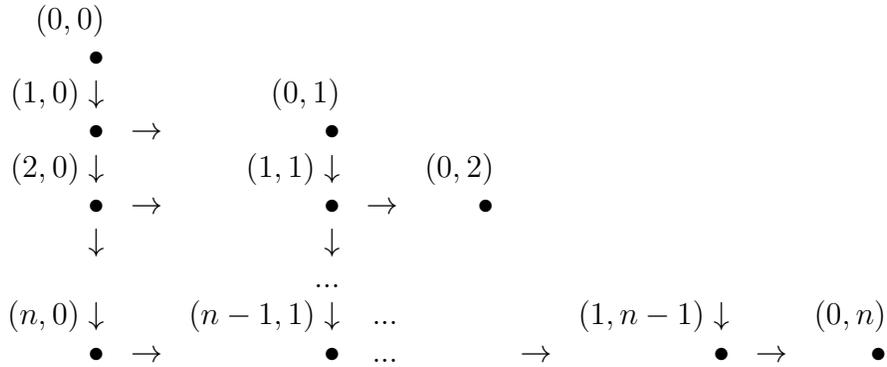
\[ \begin{array}{rrrrrrrrrrr}
  (0,0) & \\
  \bullet &  &  &  &  &  &  &  \\
  (1,0)\downarrow  &  & (0,1) &  &  &  &  &  \\
  \bullet & \rightarrow &  \bullet &  &  &  &  &  \\
 (2,0)\downarrow  &  & (1,1) \downarrow&  &  (0,2) &  &  &  \\
   \bullet & \rightarrow & \bullet & \ra & \bullet & & & & & \\
 \downarrow &   &  \downarrow & &  &  &  &  & &\\
    &  & ... &  &   &   &   &   & & \\
 (n,0)\downarrow & & (n-1,1)\downarrow & ... & & &
    (1,n-1)\downarrow &   & (0,n) \\
  \bullet & \rightarrow & \bullet & ... && \rightarrow &
                                        \bullet & \rightarrow & \bullet
 \end{array} \]

For each edge introduce edge variables: $u_{ij}$ denotes the
vertical edge variable such that whose head vertex has coordinate
$(i=\text{diagonal},j=``x\text{-axis}")$. Also $v_{ij}$ denotes
the parallel edge variable whose tail vertex has coordinate
$(i,j)$. Edges will be identified with corresponding variables.

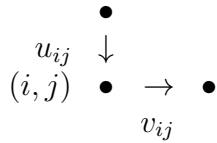
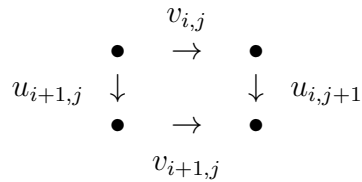
\[ \begin{array}{rrrrr}
   &\bullet &  && \\
  u_{ij} &\downarrow & && \\
    (i,j)& \bullet & \rightarrow & \bullet \\
         &      &   v_{ij}         & &
         \end{array}  \]

%

For each box

\[ \begin{array}{ccccc}
\   &   \     &    v_{i,j}   &  \      & \  \\
\   & \bullet & \rightarrow & \bullet & \  \\
 u_{i+1,j} & \downarrow & \   & \downarrow   & u_{i,j+1 } \\
\   & \bullet & \rightarrow & \bullet & \  \\
\   &  \      &     v_{i+1,j} &   \     & \
 \end{array} \]

in the graph, impose \lq\lq box\rq\rq  relation

\[ v_{i,j} u_{i,j+1} - u_{i+1,j}v_{i+1,j}=0. \]

For each roof

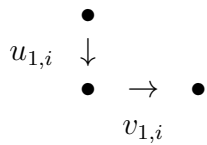
\[ \begin{array}{cccc}
    &\bullet &  & \\
u_{1,i} &\downarrow & &\\
    & \bullet & \rightarrow & \bullet \\
         &      &   v_{1,i}          &
         \end{array}  \]
in the graph, impose \lq\lq roof\rq\rq equation
\[ u_{1,i} v_{1,i} = q_i \ . \]

For nonzero complex numbers $q_i$, $i=1,...,n$, define $Y_q$ to
be the affine variety $$\Spec \CC [{\bf u},{ \bf v}]/ ( ...,
v_{i,j} u_{i,j+1} - u_{i+1,j}v_{i+1,j} , ...,u_{1,i} v_{1,i} -
q_i ,...),$$ given by all box equations and all roof equations.
  It is a complex torus
$(\CC ^\ti )^{n(n+1)/2}$. The fact follows from the
interpretation of each edge variable as the fractions of \lq\lq
voltages\rq\rq at vertices of the edge. To each vertex, assign
free vertex variables $T_{ij}$ with one condition $\sum _{i=0}^n
T_{0i}=0$. Then  it is natural to let $u_{i,j}= e^{T_{i,j} -
T_{i-1,j}}$ and $v_{i,j}= e^{T_{i-1,j+1} - T_{i,j}}$.

Now introduce a stationary phase integral
\[ \mathcal{I}_\Gamma (t_0,...,t_n)
= \int _{\Gamma \subset Y_q} e^{f_q/\hbar}\omega _q .\]
 It is a complex valued
function on $(t_0=T_{00},t_1=T_{0,1},...,t_n=T_{0n})$ with
$t_0 + ...+ t_n =0$. Here $\Gamma$ is
a descending Morse cycle of $\mathrm{Re}f_q$ with a suitable
Riemannian metric on $Y_q$.
$\Gamma$ varies covariant constantly with the Gauss-Manin
connection on the relative
homology bundle with fibers
$H_{n(n+1)/2}(Y_q, \mathrm{Re}f_q=-\infty)$. Let
\[ \omega _q = \bigwedge _{\text{all vertices not on the main diagonal}} d T_v .\]
Observe that the form $\omega _q$ can be
defined on $Y_q$ since it is translation invariant.

To outer edges $\epsilon $, that is,
edges whose targets are $(i,0)$ or sources are $(i,n-i)$,
assign weights
$\kl_\ke$ by \[ \kl _{\ke} = \left\{  \begin{array}{cc}
         - \kl_{i-1} -\frac12\sum _{j<i-1}\kl_j
& \mbox{if $\ke = v_{i,n-i}$ }\\
 \kl _{i-1} + \frac12 \sum _{j<i-1}\kl_j    & \mbox{if $\ke=u_{i,0}$}
                 \end{array}\right. .\]

{\tiny
\[ \begin{array}{rrrrrrrr}
 \bullet &  &  &  &  &  &  &  \\
 \kl _0 \downarrow  &  &  &  &  &  &  &  \\
  \bullet & \rightarrow &  \bullet &  &  &  &  &  \\
 \kl _1 + \frac12\kl_0\downarrow &   &  \downarrow & &  &  &  &  \\

   ... &  & ... &  &   &   &   &   \\
 \kl _{n-1}+\frac12\sum _{i<n-1}\kl _i\downarrow &      & \downarrow & ... &     &  \downarrow &   & \\
  \bullet & \rightarrow & \bullet & ... & \rightarrow &
                                                \bullet & \rightarrow & \bullet   \\
           & -\kl_{n-1} -\frac12\sum_{i<n-1}\kl _i & & ... & -\kl _1 -\frac12\kl_0 & & -\kl _0
 \end{array} \]}

For each edge $\ke$ which are neither at the far left nor at the bottom of the graph,
let \[ \kl_{\ke} = \left\{  \begin{array}{cc}
         - \frac12\kl _{i-1} & \mbox{if $\ke = v_{i,j}$}\\
          \frac12\kl_{i-1}   & \mbox{if $\ke=u_{i,j}$}
                             \end{array}\right. .\]

\[ \begin{array}{cccc}
                  &\bullet & &  \\
\frac12\kl_{i-1}  &\downarrow & &\\
                  & \bullet & \rightarrow                  & \bullet \\
 (i,j)            &         &   -\frac12\kl_{i-1}          &
         \end{array}  \]

{\em The assignment of weights on edges is given to satisfy the condition that the sum of outgoing
weights minus the sum of incoming weights at each vertex at diagonal level $k$ is
exactly $\lambda _{k} - \lambda _{k-1}$, if we set $\lambda _{-1} =0$.} This property will
enable us to prove the equivariant mirror theorem for flag manifolds.

Now define a phase function   $f_q$ by
\begin{eqnarray*} f_q &=&\sum _{i>0,j} (u_{i,j}+v_{i,j})
+  \sum _{i>0,j}(\kl _{u_{i,j}} \ln u_{i,j} + \kl _{v_{i,j}} \ln v_{i,j}) \\
& =& \sum _{i>0,j}(e^{T_{i,j}-T_{i-1,j}}+e^{T_{i-1,j+1}-T_{i,j}}) \\
   && \quad + \sum _{i>0,j}
(\kl _{u_{i,j}}(T_{i,j}-T_{i-1,j}) +\kl _{v_{i,j}}(
T_{i-1,j+1}-T_{i,j})) .
\end{eqnarray*}
The first and second summation terms are respectively $f_q'$ and the weight factor
that we mentioned in the beginning of this section.

Define $\sigma _i$ by equation $x^{n+1}+\sigma _1 x^n + ... + \sigma _{n+1}
= \prod _{i=0}^n(x-\kl _i )$.
We state the equivariant mirror theorem for flag
manifolds and shall prove it in section \ref{pf1} and \ref{pf2}.

\begin{Thm}\label{sol} Let $\Gamma$ be a descending Morse cycle of $\mathrm{Re}f_q$ at
a nondegenerate critical point of $f_q$. Then the stationary
phase integral
\[ \mathcal{I}_{\Gamma} (t_0,...,t_n)= \int _{\Gamma \subset Y_q} e^{f_q/\hbar} \omega _q \]
satisfies eigenvalue differential equations $D_1\mathcal{I} =
\sigma _1 \mathcal{I}$, $D_2\mathcal{I} = \sigma _2\mathcal{I}$ ,
..., $D_{n+1}\mathcal{I} = \sigma _{n+1} \mathcal{I}$.
\end{Thm}

\begin{Cor}\label{sol2} Let $q$ be a general point in $\CC ^{n}$. Then there are
$(n+1)!$ critical points of $f_q$ . Furthermore they are all
nondegenerate.
\end{Cor}

\bigskip
{\em Remark.} When $\kl _i =0$ for all $i$, the above theorem is proven
in \cite{Gs}.

\subsection{Proof of theorem \ref{sol}}\label{pf1}

Since $\frac{\partial f_q}{\partial t_i}$
does not depend on $t_j$ if $j\ne i$,
the amplitude created by operating $D_i$ on $\mathcal{I}_\Gamma$
are the corresponding coefficients
of the characteristic polynomial of $A_{1} - \kl _0 I$, where
 \[ A_{k} =\left[\begin{array}{cccccc}
-u_{k0}  & u_{k0}v_{k0} & 0 & ...  & & \\
-1 & v_{k0}-u_{k1}  & u_{k1}v_{k1} & ... & & \\
 & & ... & & & \\
  & &  & & & u_{k,n-k}v_{k,n-k} \\
  & & & & -1 & v_{k,n-k}
  \end{array}             \right]
  \]
for $k=1,...,n$. Hence
\begin{eqnarray*}
&& (x^{n+1} + D_1 x^n + ...+ D_{n+1} ) \cdot \int _{\Gamma} e^{f_q/\hbar} \omega _q
\\
&=& \int _{\Gamma} \det ( A_1 - \kl _0 I +xI)
\cdot  e^{f_q/\hbar}\omega _q .
\end{eqnarray*}
To prove the theorem it is enough to show  the following.

\begin{Prop}\label{UV} For $k=1,...,n$
\[  \int _{\Gamma}\det ( A_{k} - \kl _{k-1} I +xI)\cdot  e^{f_q/\hbar}\omega _q
 = (x-\kl _{k-1})
\int _{\Gamma}\det (A_{k+1} - \kl _kI+xI) \cdot e^{f_q/\hbar }\omega _q,\]
where $A_{n+1}=0$.
\end{Prop}
{\em Proof.}
As noticed in \cite{Gs}, $A_k=U_kV_k$ where
\[U_k=\left[\begin{array}{ccccccc}
u_{k0}  & 0    &   & ...  & & \\
1    & u_{k1}  & 0 & ... & & \\
  &  1   & u_{k2}  & ... & & \\
     &      & ... & & & \\
     &      &  &      &1 & u_{k,n-k} & 0\\
     &      &  &       &  & 1 & 0
  \end{array}             \right]  \]
and \[
  V_k =\left[\begin{array}{cccccc}
-1  & v_{k0} & 0 & ...  & & \\
0 & -1  & v_{k1} & ... & & \\
0 & 0 & -1 & ... & & \\
 & &     ... & & & \\
  & &  & & & v_{k,n-k} \\
  & & & & 0 & -1
  \end{array}             \right]     .
  \]
Since $V_k$ is invertible, the characteristic polynomial of
$A_{k}-\kl _{k-1} I$ is equal to that of $V_kU_k-\kl _{k-1}I$.  Using  the
critical point conditions $$\frac{\partial
f_q}{\partial T_{k,i}}
=u_{ki}-v_{ki}+v_{k+1,i-1}-u_{k+1,i}+\kl _{k-1}-\kl _k$$
and box relations
$u_{k,i+1}v_{ki}=u_{k+1,i}v_{k+1,i}$, we see that
 \[V_kU_k-\kl _{k-1} I = \left[\begin{array}{cc}
 A_{k+1} - \kl _k I & 0 \\
 0,...,0,-1 & -\kl_{k-1}
 \end{array}\right] - \mathrm{diag} (\frac{\partial f_q}{\partial T_{k,0}},...,
 \frac{\partial f_q}{\partial T_{k,n-k}}, 0). \]

Note that the terms in $\det (V_kU_k-\kl_{k-1})$ which involve
the derivatives of $\frac{\partial f_q}{\partial T_{k,i}}$ might
have extra multiples of  edge variables whose vertices are not
$(k,i)$. However if $P$ is a polynomial of edge variables except
whose vertices are not $(k,i)$, the following holds:
\begin{eqnarray*}\int _{\Gamma} \frac{\partial f_q/\hbar}{\partial T_{k,i}}
P e^{f_q/\hbar} \omega _q &=&\int _{\Gamma} \frac{\partial
}{\partial T_{k,i}}( P e^{f_q/\hbar})\omega_q\\
&=& \int _{\Gamma} L_{\frac{\partial}{\partial T_{k,i}}}(
P e^{f_q/\hbar} \omega _q)\\
&=&\int_{\Gamma} d i_{\frac{\partial}{\partial T_{k,i}}} P
e^{f_q/\hbar}
\omega _q\\
&=& 0.\end{eqnarray*} This completes the proof.

\subsection{Proof of Corollary \ref{sol2}}\label{pf2}

Let $$ Y = \{ (...,u_{ij},v_{i,j},...)| \text{all } u_{i,j} \ne
0, \text{all } v_{i,j}\ne 0, \text{all box relations}\}.$$ $Y$ is
a complex torus $(\CC ^\ti)^{n(n+1)/2+n}$. For $\kl=(\kl
_0,...,\kl_n )$ with $\sum \kl _i =0$, denote by $Z_{\kl}$ the
closed subscheme in $Y$ defined by the ideal generated by
$\frac{\partial f_q}{\partial T _{ij}}$, $\forall \ i > 0$, $j$.
Let $X_{\kl} = \Spec\CC [p_0,...,p_n, q_1^{\pm 1},...,q_n^{\pm 1}
] / (D_1 (p,q)-\sigma _1 ,...,D_{n+1}(p,q)-\sigma _{n+1})$, which
is an irreducible nonsingular rational variety \cite{Ko}. Since
$\det (A_{1}-\kl _0 I +xI)= \Pi _{i=0}^n (x-\kl _i)$ on $Z_\kl$ as
shown in the proof of proposition \ref{UV}, a morphism $\phi:
Z_\kl \ra X_\kl$ is defined by $p_0 = -u_{1,0}$, $p_1 = v_{1,0}
-u_{1,1}$, ..., $p_n=v_{1,n-1}$,
$q_1=u_{1,0}v_{1,0}$,...,$q_n=u_{1,n-1}v_{1,n-1}$. Direct
inverting of $\phi$ shows that $\phi$ is a birational morphism. Let
$\pi$ be the projection from $X_\kl$ to $(\CC ^\ti)^n$ defined by
$\pi (p,q)=q$. The projection is a finite map of degree $(n+1)!$
which is verified by the dimension of the algebra $\CC
[p,q]/D_1(p,q)-\sigma _1,...,D_{n+1}(p,q)-\sigma_{n+1})$ at $q=0$.
Now the proof of the corollary follows from the irreducibility of
$X_\kl$, the degree of the projection $\pi$ and Sard's lemma.

\subsection{Study of the phase function $f_q$}

Given a sequence $\sigma = (k_1,...,k_n)$ of integers
where   $0\le k_i \le n-i+1$,
let $w^\sigma _{ij}= u_{ij}$ if $j<k_i$, otherwise
let $w^{\sigma} _{ij} = v_{ij}$.
We shall express $f_q$ in terms of these independent variables $w^{\sigma}_{ij}$.

We associate a sequence $\rho ^{\sigma}_{ij}$ of weights  to
$\sigma$, inductively on $n-i$ and then $j$, where $i=1,...,n+1$ and
$j=-1, 0,..., n-i+1$. First, let $\rho ^\sigma _{n,-1}=0$,
$\rho ^\sigma _{i,n-i+1}=\sum _{k=i-1}^n \kl _k$. And let
$\rho ^{\sigma} _{ij}=\rho ^{\sigma} _{i+1,j}$
if $w^\sigma_{ij}=u_{ij}$,
otherwise let $\rho ^{\sigma} _{ij} = \rho ^{\sigma} _{i+1,j-1} +\kl_{i-1}$.
It is easy to see inductively that
$\{ \rho ^{\sigma} _{ij} -
 \rho ^{\sigma} _{i,j-1} \ | \ j=0,1,..,n-i+1\} = \{ \kl _{i-1},...,
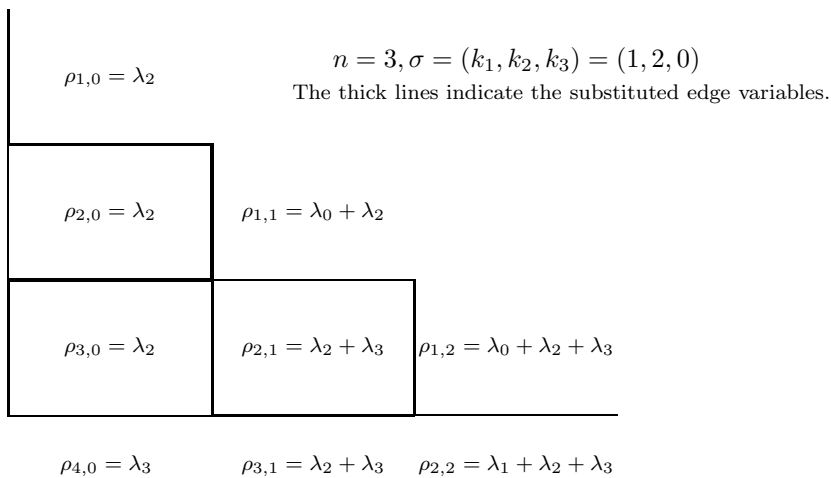
\kl _n\}$.

\setlength{\unitlength}{0.9cm}
 \begin{picture}(9,8)
 \thinlines
 \put(1,1){\framebox (3,2){\tiny$\rho_{3,0}=\lambda_2$}}
 \put(4,1){\makebox (3,2){\tiny$\rho_{2,1}=\lambda_2+\lambda_3$}}
 \put(1,3){\framebox (3,2){\tiny$\rho_{2,0}=\lambda_2$}}
 \put(4,3){\line(1,0) {3}}
 \put(7,1){\line(1,0){3}}
 \put(1,5){\line(0,1){2}}
 \put(1,5){\makebox(3,2){\tiny$\rho_{1,0}=\lambda_2$}}
 \put(4,3){\makebox(3,2){\tiny$\rho_{1,1}=\lambda_0+\lambda_2$}}
 \put(7,1){\makebox(3,2){\tiny$\rho_{1,2}=\lambda_0+\lambda_2+\lambda_3$}}
 \linethickness{5mm}
\thicklines \put(1,1){\line(0,1){2}} \put(1,3){\line(1,0){3}}
\thicklines \put(4,3){\line(0,1){2}} \put(1,5){\line(1,0){3}}
\thicklines \put(4,1){\line(1,0){3}}\put(7,1){\line(0,1){2}}
\put(1,-0.5){\makebox(3,1.5){\tiny$\rho_{4,0}=\lambda_3$
}}\put(4,-0.5){\makebox(3,1.5){\tiny$\rho_{3,1}=\lambda_2+\lambda_3$}}
\put(7,-0.5){\makebox(3,1.5){\tiny$\rho_{2,2}=\lambda_1+\lambda_2+\lambda_3$}}
\put(5, 6){\makebox(7,0.5) {\footnotesize  $ n=3, \sigma= (k_1,k_2,k_3)=(1,2,0)$}  }

\put(5.5, 5.5){\makebox(7,0.5){\tiny  \quad The thick lines indicate the substituted edge variables.}}
 \end{picture}

\begin{Prop}\label{phaseft}
The phase function $f_q$ is, in terms of $w^\sigma _{ij}$,
\begin{eqnarray*}
& & \sum _{i=1,...,n} \rho ^{\sigma}_{1,i-1}\ln q_i + \sum _{i=1,...,n,
j=0,...,n-i} ( w^\sigma _{ij} +  r_{ij}
 (... ,w^\sigma _{ab},...,q_1,...,q_n )
\\
&+& \sigma (i,j)
\ln w^\sigma _{ij} ),\end{eqnarray*}
where $\sigma (i,j)$ is
$\kl _{i-1} - (\rho ^{\sigma}_{i,j} - \rho ^{\sigma} _{i, j-1}) $
if $j < k_i$ and $\sigma (i,j)$ is
$-\kl _{i-1} + (\rho ^{\sigma}_{i,j+1} - \rho ^{\sigma} _{i, j}) $
otherwise. Here $r_{ij}$ are monomials in $(w^\sigma _{ab})^{\pm 1}$
and $q_k$ with at least one factor among $q_1,...,q_n$.
 \end{Prop}

 {\em Proof.\/}
To apply induction on $n$, let $L$ be the triangular graph
introduced in the beginning of section 3 and
let $L_k$ be the triangular sublattice of $L$
whose vertices $(i, j)$ satisfy
inequalities $k\le i $. The edges of $L_k$ are by definition
the edges of $L$, connecting neighboring vertices of $L_k$.
Note that $L_0$ is the full lattice $L$. Introduce
new variables associated to  the sublattice $L_k$ as follows, $$
Q^k_{j+1}=u_{k+1,j}v_{k+1,j} \quad \, \Delta _k =\sum_{j=0}^{k-1}
\lambda_j .$$ Moreover let $ Q^0_j=q_j , \Delta_0=0$. Consider
$$(f_k)_{Q^k} = \sum_{\mathrm{vertices} \ \in L_k} (u_{i,j}+v_{i,j})
+\sum_{\mathrm{vertices} \ \in L_k} (\kl _{u_{i,j}}\ln u_{i,j} + \kl _{v_{i,j}}\ln v_{i,j}).$$
By the induction hypothesis, we can  assume that
\begin{eqnarray*} &(f_k)_{Q^k} = & \sum _{i=1,...,n-k} (\rho ^{\sigma}_{k+1,i-1}+\frac{\Delta_k}{2})\ln
Q^k_i \\& +& \sum _{i=1,...,n-k ;\ j=0,...,n-i} ( w_{ij} +
r^k_{ij} (... ,w^\sigma _{ab},...,Q^k_1,...,Q^k_{n-k} )
\\
&+& \sigma (i,j)
\ln w^\sigma _{ij} ) ,\end{eqnarray*} where
$r^k_{ij}$ are monomials in $(w^\sigma _{ab})^{\pm 1}$ and
$Q^k_1$, ..., $Q^k_{n-k}$. Then substitute  $ v_{k,i} u_{k,i+1}$
for $Q^k_{i+1}$.  Then derive
$(f_{k-1})_{Q^{k-1}}$ in terms of $w^\sigma_{i,j}$ and $Q^{k-1}_j$. The rest is straightforward. $\Box$

 \bigskip

Let $f^\sigma_q = f_q - \sum _i \rho ^\sigma _{1,i-1} \ln q_i$.
Then $f_q^\sigma$ is regular at $q_1=...=q_n=0$ and has exactly one
simple critical point at the origin $q=0$.
Therefore there is exactly one critical point of $f_q$ whose limit as $q$ goes to 0
is the simple critical point of $f^\sigma _0$. Let us denote by $\Gamma _\sigma$
 the descending Morse
cycle  of $f_q$, whose limit is the Morse cycle $\Gamma _{\sigma ,0}$
of $f_0^\sigma$. Now we obtain  a corollary which will be applied to
the proof of theorem \ref{Rconjforflag}.

\begin{Cor}\label{asy}
   \[ \lim _{q\ra 0}e^{-\sum _{j=0}^n
(\rho ^{\sigma} _{1,j-1}-\rho ^{\sigma} _{1j})t_j }
 \mathcal{I}_{\Gamma _\sigma}
= \prod _{i=1,...,n; \ j=0,...,n-i }
\int _{C_{\sigma(i,j)}} e^{w_{ij}/\hbar }
 w_{ij}^{\sigma (i,j) /\hbar }\frac{dw_{ij}}{w_{ij}},\]
 where $\sigma (i,j)$ is defined as in the previous proposition
 and $C_{\sigma(i,j)}$ is the descending Morse cycle of an one-variable
 function
$x+ \sigma (i,j) \ln x$.
\end{Cor}
\bigskip

For later use, for a given sequence $\sigma = (k_1,...,k_n)$ of integers $0\le k_i \le n-i+1$,
we define $\sigma (i)\in \{0,1,...,n\}$ by the requirement
$$\kl _{\sigma(j)}=\rho ^{\sigma} _{1,j}-\rho ^{\sigma} _{1,j-1} ,$$ so
that we may identify $\sigma$ as a permutation element of the symmetric group $S_{n+1}$ of $n+1$ letters.

\section{the Virasoro Conjecture}

\subsection{The fundamental solution of the form $\Psi R e^{U/\hbar}$}

Now we  recall the equivariant big quantum cohomology
of a projective algebraic manifold $X$ with a Hamiltonian torus $T$-action.
Fix a homogeneous basis $\{\phi _\ka \}$ of $H^*_T(X, \CC)$.
For simplicity assume that $H^*(X,\CC ) = H^{\mathrm{even}}(X,\CC)$ (i.e. there
are no odd classes) and there is a basis $\omega _i$ of $H^2(X,\ZZ)$ such that
$<\omega _i, [C]>$ is nonnegative for any curve (effective) class $[C]$.
Introduce formal parameters $Q_i$ and for $d\in H_2(X,\ZZ)$
let $Q^d = \prod Q_i^{d_i}$ where $d_i=<\omega _i,d>$.
A potential function $F(t)$ is by definition
\begin{eqnarray*} F(t) &= &\sum _{m,d} \frac{Q^d}{m!}
\int _{[\overline{M}_{0,m}(X, d)]^{\mathrm{virt}}} ev_1^*(t ) ...
ev_m^* (t).
\end{eqnarray*}
Here is an explanation of notation in the equation. First,
$[\overline{M}_{g,m}(X, d)]^{\mathrm{virt}}$ denotes the virtual
fundamental class of the moduli space $\overline{M}_{g,m}(X, d)$
of stable maps to $X$ with genus $g$, $m$ marked points and
degree $d\in H_2(X,\ZZ )$. Secondly, $t\in H^*_T(X)$ is considered
as $t=\sum t_\ka\phi _\ka $ and the integral is taken as the
equivariant pushforward so that the value is in $H^*(BT)$. The
potential function defines the equivariant big quantum cohomology
of $X$ as following. Let $H= H^*_T(X,\CC [[Q]])$. The quantum
cohomology is a certain multiplication structure on the tangent
space $T_tH$. The quantum product $\frac{\partial}{\partial
t_\ka} \circ _{t}\frac{\partial}{\partial t_\kb}$ is defined by
equation
$$( \frac{\partial}{\partial t_{\ka}} \circ _{t} \frac{\partial}{\partial t_{\kb}} ,\frac{\partial}{\partial t_{\kg}})
= \frac{\partial}{\partial t_{\ka}}\frac{\partial}{\partial
t_{\kb}}\frac{\partial}{ \partial t_{\kg}} F (t).$$ Here $(,)$ is
the equivariant Poincar\'e pairing on $T_tH$ defined by
$(\frac{\partial}{\partial t_{\ka}}, \frac{\partial}{\partial
t_{\kb}}):=\int _X \phi _\ka \wedge \phi _\kb =:g_{\ka\kb}$.
There is a pencil of flat connections $\nabla _{\hbar} = \hbar d
- \sum _{\ka} dt_\ka\frac{\partial}{\partial t_\ka} \circ$ with
parameter $\hbar$. Here $d$ is the Levi-Civita connection with
respect to the vector space structure of $H$. (Strictly speaking,
$\nabla _{\hbar} /\hbar$ are  connections.)

$\Spec (Vect (H))$ defines a Lagrangian formal scheme
$L\subset T^*H$ over the ground ring $\CC [\kl _0,...,\kl _n][[Q]]$, where
$H^*(BT)=\CC [\kl _0,...,\kl_n]$. Let $\mathrm{rank} H^*(X) = N$ and
let $u$ be a semi-simple point in $H$,
that is, $T_u H$ is a semi-simple algebra.
Then $L_u \in T^*_u H$
consists of $N$ many simple points and there are $N$ many
local functions $u_{1},...,u_{N}$
at $u$ such that the canonical 1-form
$\sum _{\ka} P_{\ka} dt_{\ka}$
restricted to $L_{i}$ coincides with
$du_i :=\sum \frac{\partial u_i}{\partial t_\ka}dt_\ka$,
where $L_{i}$ are branches of $L$. Since $\{ du_i \}$ is linearly independent,
 $\{ u_{i} \}$ forms a local coordinate
system of $H$ and the dual vector fields $\frac{\partial}{\partial u_i}$ can be defined.
 Notice that $\frac{\partial}{\partial u_{i}} \circ _{t} \frac{\partial }{\partial u_{j}}
 = \delta _{ij}\frac{\partial }{\partial u_{j}} $.
We call $u_{i}$ canonical coordinates.
These coordinates are unique up to order and
additions by constants.

Let $U$ be $\mathrm{diag} (u_{1},...,u_{N})$ and $A^1$ be the
connection 1-form matrix of $\nabla _{\hbar}$ so that $\nabla
_{\hbar} \frac{\partial}{\partial t_\ka} = \sum _{\kb} A^1_{\kb,
\ka } \frac{\partial}{\partial t_\kb}$.  Let $1/ \Delta _\kb$ be
the length of $\frac{\partial}{\partial u_{\kb }}$ with respect to
the Poincare paring metric and let $\Delta _\ka \frac{\partial
}{\partial u_{\ka }} = \sum _\kb \Psi _{\kb , \ka} \frac{\partial
}{\partial t_\beta }$. Then the computation of connection 1-form
$\hbar\Psi ^{-1} d\Psi -dU$ with respect to frame $\Delta _{\ka}
\frac{\partial}{\partial u_{\ka}}$ using the expression $\nabla
_{\hbar} = \hbar d - \sum du_{\kb}\frac{\partial}{\partial u_\kb
}\circ$ shows that $$ A^1 \Psi =\Psi dU.$$  Notice that $u_i$
linearly depends on $t_\ka$'s modulo $Q$ and $\Psi$ is constant
modulo $Q$.

\bigskip

{\em Example.} Let $X=Fl(2)=\CC P^1$. Consider $T=\CC ^{\ti}$
action as $t\cdot [x_0,x_1]=[tx_0,t^{-1}x_1]$ where $t\in T$ and
$[x_0,x_1]\in \CC P^1$. The Lagrangian scheme $L\subset T^*H$ is
defined by equation $p^2 = \kl ^2 + e^tQ $ and ${\bf 1}=1$ where
$p$ is the cotangent direction variable associated to $-t_0=t_1$
and ${\bf 1}$ is the cotangent direction variable associated to
the identity element of $H^*_T(\CC P^1)$. The fiber $L _{(t_0,t)}
$ of $L$ over a point $(t_0,t)$ of $H$ consists of two points $\{
(1, \pm \sqrt{\kl ^2+e^{t}Q}) \}$. Then $du_{\pm} = dt_0 +
p_{\pm} dt =dt_0\pm \sqrt{\kl ^2 +e^tQ} dt$ so that we may set
$u_\pm =t_0 \pm \int \sqrt{\kl ^2 +e^tQ} dt$. The canonical
vector fields are $\frac{\partial }{\partial
u_{\pm}}=\frac{1}{2p_{\pm}}(p_{\pm}\frac{\partial }{\partial t_0}+
\frac{\partial}{\partial t})$. Let $v=\sqrt{p_+}$. Then with
respect to the ordered basis $\{ 1, p\}$ of $H^*(\CC P^1)$

\[ A^1 =
\left[\begin{array}{cc}
 dt_0  &  v^4 dt\\
 dt &  dt_0
 \end{array}\right] , \
 \Psi =  \frac{1}{\sqrt{2}}
\left[\begin{array}{cc}
v & -iv \\
1/v & i/v
\end{array}\right]. \]

On the other hand, we have shown that $\{ (p,q)| \frac{\partial
f_q}{\partial y_\ke } =0, p= \frac{\partial f_q}{\partial t}\}$
is the Lagrangian scheme $L'=\Spec \CC [ p ,q ^{\pm 1}] / (p^2 -
q -\kl ^2 )$. At critical points of $f_q$, $df_q = p dt$, so that
$u_\pm (\kl , q) = t_0 + f_q(p_\pm (\kl ,q),q)$ (up to constant
addition). Thus, $e^{t_0/\hbar}\mathcal{I}_{\Gamma _{\pm}} \sim
e^{u_\pm /\hbar}$ (up to constant multiplication). In fact a
direct computation shows that at critical points
$$f_q (p_\pm,q) = 2p_\pm + \kl _0 \ln (\kl _1 + p_\pm ) + \kl _1 \ln (\kl _0 + p _\pm)$$
and $\frac{\partial f_q(p_\pm,q)}{\partial t} = p_\pm$. $\Box$
\bigskip

Near a semi-simple point $u$, due to the papers \cite{Gelliptic, Gsemi, Gh} there is a
fundamental solution
of form $$\Psi R e^{U/\hbar}$$
such that $\Psi$ and $U$ are defined as above and
$$R(\hbar) = 1+ R_1 \hbar +R_2 \hbar ^2 + ...,$$
with $R^t(\hbar) R(-\hbar ) =1$. Since $Q_i
\frac{\partial}{\partial Q_i} F (t)=q_i\frac{\partial}{\partial
q_i} F(t)$,  we require that $Q_i \frac{\partial}{\partial Q_i}
R=q_i\frac{\partial}{\partial q_i} R$, where $e^{t_i}=q_i$.
Such matrix $R$ exists
uniquely up to right multiplication by diagonal matrices $\exp
(\sum _{k=1}^{\infty} a_k\hbar ^{2k-1} )$ where $a_k
=\mathrm{diag} (a_k^{(1)},...,a_k^{(N)})$ are constants.

Suppose $T$ acts on $X$ with isolated fixed points and  let $w_i$,
$i=1,...,N$ be $T$- fixed points of $X$ and let $N_l^{(i)}=
\sum_j \frac{1}{\chi _j (w_i)^l}$ where $\chi _j(w_i)$  are the
weights of the induced torus action on the cotangent space
$T^*_{w_i} X$ at $w_i$.

\bigskip
Let $u$ be a semi-simple point of $H^*_T(X)$. There is an unique
asymptotic fundamental solution of type $\Psi R e^{U/\hbar}$ of
the connection $\nabla _{\hbar}$ satisfying the following four
conditions on $R=1+R_1\hbar + R_2\hbar ^2 + ...$:

0. the divisor condition: $Q_i\frac{\partial}{\partial Q_i} R =
q_i\frac{\partial}{\partial q_i}R$ for all $i$,

1. the orthogonal condition: $R^t(\hbar)R(-\hbar)=1$,

2. the classical limit condition: its classical limit of letting
$Q\ra 0$ is $\exp (\mathrm{diag} (b_1,...,b_N))$ where
\[ b_i(\hbar ) = \sum _{k=1}^{\infty} N_{2k-1}^{(i)}
\frac{B_{2k}}{2k} \frac{\hbar ^{2k-1}}{2k-1}.\] Here $B_{2k}$ are
Bernoulli numbers defined by $x/(1-e^{-x}) = 1+x/2+ \sum
_{k=1}^{\infty} B_{2k}x^{2k}/(2k)!$,

3. the equivariant homogeneity condition: $(\hbar \partial
_{\hbar} + \sum u_{\ka}
\partial _{\ka}+ \sum\lambda_i\partial_{\lambda_i})R(\hbar) =0$.
\bigskip

{\bf R-Conjecture} {\rm (\cite{Gh})} The matrix $R$ has the
nonequivariant limit of letting $\kl _0\ra 0,..., \kl_n\ra 0$ and
thus the nonequivariant limit satisfies the homogeneity condition
$(\hbar
\partial _{\hbar} + \sum u_{\ka} \partial _{\ka} )R(\hbar) =0$
where $u_\ka$ denote the corresponding nonequivariant canonical
coordinates.

\bigskip


{\em We shall prove the $R$-conjecture restricted to $H^2_T(X)$}
when $X$ is a flag manifold. In such restriction the homogeneity
condition is replaced by $(\hbar\frac{\partial}{\partial
\hbar}+\sum_i < c_1(T_X), d>\frac{\partial }{\partial t_i} )R
(\hbar)|_{H^2_T(X)} = 0$.
 There are two known ways to obtain fundamental
solutions of form $\Psi R e^{U/\hbar}$ restricted to  $
H^2_T(X)$. One is the localization expansion of the construction
by two-pointed gravitational Gromov-Witten invariants and the
other one is the stationary phase approximation of the mirror
construction \cite{Gelliptic}. We use both. There is a natural 1-1
correspondence between the fixed points $w_\ka$ and canonical
coordinates $u_\ka$ since the points can be identified with
points $L_t$ of the Lagrangian variety $L$
 (first with $Q=0$ and general $\kl$ and then continuously).
 So we use the same label.
If $$S_{\kb , \ka} = \sum_d \frac{Q^d}{m!}
\int_{[\overline{M}_{0,m+2}(Fl(n+1),d)]^{\mathrm{virt}}}
ev_1^*(\phi _\kb)(\prod _{i=2}^{m+1} ev_i ^* (t))
\frac{ev_{m+2}^*(\phi_\ka)}{\hbar -\psi _{m+2}}$$ then
 $$s_{\ka}=\sum_{\kb, \kg }S_{\kb, \ka} g^{\kb \kg} \phi _\kg$$  is a flat
section for all $\ka$ \cite{Ge}. Recall if $\phi _{\kb} =1$,
$t=\sum t_i\omega_i \in H^2_T(X)$ then $$S_{\kb ,\ka} = \sum
_{d} Q^d q^d \int _{[\overline{M}_{0,1}(X,d)]^{\mathrm{virt}} }
\frac{ev_1^*(\phi_\ka \wedge \exp (t/\hbar)) ) }{\hbar(\hbar
-\psi _1)}, $$ which is  denoted by $S_{1,\ka}$. We  see the
equality using the string and divisor axioms.

Suppose that {\em $X$ is a Fano manifold}. Let $T$ act on $X$ and
{\em assume that $H^*_T(X)$ is
generated by the second cohomology classes} $\{ \omega _i \}$. Here $T$ could be
just the trivial action. Let $\mathcal{D}$ be the algebra of
differential operators $\hbar\frac{\partial}{\partial t_\ka}$ over
$\CC [\kl, Qq, \hbar]$.   Then the following lemma holds.
\begin{Lemma}
For each $\phi \in H^*(X)$, there is a differential operator $D_{\phi}$ in $\mathcal{D}$ such that
$(s_\ka , \phi ) = D_{\phi} S_{1,\ka} (t)$ for every $\ka$, $t\in H^2_T(X)$.
\end{Lemma}
{\em Proof.} We may set $D_{\omega _i}= \hbar\frac{\partial}{\partial t_i}$.
If $D_{\phi } S_{1,\alpha} = (s_\ka, \phi )$
and $D_{\phi  \circ \omega _i - \phi \wedge \omega _i} S_{1,\ka } = (s_\ka,
    \phi \circ \omega _i - \phi \wedge \omega _i)$, then
  $(\hbar\frac{\partial}{\partial t_i} D_{\phi } - D_{\phi  \circ \omega _i - \phi \wedge \omega _i})
 S_{1,\ka} = (s_\ka, \phi\wedge \omega _i)$ by the flatness of the section $s_\ka$,
 where the vector fields on the affine space $H^*(X)$ are identified with classes in $H^*(X)$.
Notice that the degree of each homogeneous term of the $q$-linear combination class
$\phi \circ \omega _i - \phi \wedge \omega _i$
is less than $\deg  \phi \wedge \omega _i$ since $X$ is a Fano manifold. This completes the proof.
$\Box$

\bigskip

{\em Remark.}
When $X=Fl(n+1)$, since the oscillatory integrals $\mathcal{I}$ and $S_{1,\kb}$ with $Q=1$
generated the same $D$-module, we conclude that the matrix $(D_{\ka} \mathcal{I}_{j})$ is a fundamental solution of
$\nabla _{\hbar}$ with $Q=1$ for some $D_\ka \in \mathcal{D}$.

\begin{Thm}\label{Rconjforflag} The $R$-conjecture restricted to $H^2_T(X)$ is true for $X=Fl(n+1)$ with the $T$-action.
\end{Thm}
{\em Proof.\/}
We may let $Q=1$ by the requirement 0 of the conjecture.
Recall that $s_i=\sum _{\ka , \kb}(D_\ka \mathcal{I}_{\Gamma _i}) g^{\ka , \kb}\frac{\partial}{\partial t_\kb}$
is a flat section by remark above. Since
$$D_\ka \mathcal{I}_{\Gamma _i} = \int _{\Gamma _i} e^{f_q/\hbar}\phi _{\ka, q} \omega $$ for some polynomial
$\phi _{\ka , q}$ in $\kl$, $\hbar$, and $q_i$ (see Proposition
\ref{phaseft}), we have the asymptotic expansion
$$e^{u_i/\hbar }
\frac{\hbar ^{(n+1)!/2} \phi _{\ka,q}(\mathrm{crit} _\sigma)}{\sqrt{\det\mathrm{Hessian}
f_q (\mathrm{crit}_{\sigma})}}
(1 + o(\hbar ))$$ of the integral. From the expansion we obtain the fundamental solution
$\Psi _{\mathrm{osc}}  R_{\mathrm{osc}} e^{U/\hbar}$
where $R_{\mathrm{osc}} = 1 + (R_1)_{\mathrm{osc}}\hbar+...$
and $$(\Psi _{\mathrm{osc}})_{\ka , \sigma} =\frac{\phi _{\ka,q}(\mathrm{crit} _\sigma)}{\sqrt{\det\mathrm{Hessian}
 f_q (\mathrm{crit}_{\sigma})}}. $$
 We shall see that $\Psi _{\mathrm{osc}}$  coincides with $\Psi$ up to
 left multiplication by a constant matrix (not depending on $\hbar$ and $q$). Also
we shall show that $R_{\mathrm{osc}}$ is the $R$ satisfying all the properties in the theorem/conjecture.

First we prove that $\Psi _{\mathrm{osc}}$ coincides with $\Psi$ modulo right multiplication
of a constant diagonal matrix, using the symmetry of the differential equation.
Since $\Psi _{\mathrm{osc}}R_{\mathrm{osc}}e^{U/\hbar}$ is
a fundamental solution, $\Psi _{\mathrm{osc}}$ is
an eigenvector matrix of the connection matrix $A^1$.
However there is a symmetry of the differential equation:
 \begin{eqnarray*} && \hbar d  (s_i (-\hbar), s_j (\hbar))
 = (\hbar d s_i(-\hbar), s_j(\hbar)) +(s_i(-\hbar),\hbar d s_j(\hbar))              \\
&=& -(A^1\Psi_{\mathrm{osc}}R_{\mathrm{osc}} (-\hbar) e^{-U/\hbar})^t G\Psi_{\mathrm{osc}} R_{\mathrm{osc}}(\hbar) e^{U/\hbar}    \\
 && + (\Psi_{\mathrm{osc}}R_{\mathrm{osc}} (-\hbar) e^{-U/\hbar})^t G A^1\Psi_{\mathrm{osc}}R_{\mathrm{osc}}(\hbar) e^{U/\hbar}
  =  0 \end{eqnarray*} since $GA$ is symmetric, where $G=(g_{\ka,\kb})$.
This shows that $\Psi_{\mathrm{osc}} ^t G\Psi_{\mathrm{osc}} $ is constant in $q$ and so the claim is true.
Now considering
$c_iD_\ka \mathcal{I}_{\Gamma _i}$ for some constant $c_i$ not depending on $\hbar$ and $q$,
we conclude that $\Psi R_{\mathrm{osc}}e^{U/\hbar}$ is a fundamental solution matrix.

We investigate $R_{\mathrm{osc}}$ modulo $q$.
In order to select $R_{\mathrm{osc}}$ not $\Psi  R_{\mathrm{osc}}$
from the phase integrals
we remove $\Psi$
by expressing the fundamental solution with respect to the basis
$\{ \Delta _\ka \frac{\partial }{\partial u_{\ka}}\}$.
If $s$ is a flat section of $\nabla _{\hbar}$, then
$(\hbar \frac{\partial}{\partial t_i} )^j (s , 1)
= ( s , p_i^j )$ modulo $q$. So, if $s= \sum b_i \Delta _i \frac{\partial}{\partial u_i}$, then
$b_\sigma = \hat{L}_\sigma (s,1)$ modulo $q$
where  $$\hat{L}_\sigma
= \prod _{i>j} \frac{\hbar\frac{\partial}{\partial t_{\sigma (i)}} - \kl _j}{\sqrt{\kl_i - \kl _j}}.$$
(Here $\sigma$ as an element of permutation is defined in the remark below of corollary
\ref{asy}.)
Therefore, the entries of $R_{\mathrm{osc}}$ that we obtained above satisfy
$$ (\hat{L}_{\tau}c_{\sigma} \int _{\Gamma _\sigma} e^f \omega) _{\tau , \sigma } \sim
R_{\mathrm{osc}}e^{U/\hbar }  \ \ \ \  (\mathrm{mod}\ q).$$ By
corollary \ref{asy}
  \[\hat{L}_{\tau}c_{\sigma}\mathcal{I}_{\Gamma _\sigma } |_{q=0}
\sim  \delta _{\sigma, \tau} e^{-\sum  \kl _{\sigma (i)} t_i/\hbar}
\prod _{n\ge i>j\ge 0} \Gamma (\frac{1}{\hbar}(\kl _{\sigma (i)} - \kl _{\sigma(j)} )) \]
up to an irrelevant multiplication factor.
However, for large $z$, \[\ln\Gamma (z) \sim
(z-\frac12) \ln z  - z + \frac12 \ln 2\pi + \sum _{i=1} ^{\infty}
\frac{B_{2i}}{2i(2i-1) z^{2i-1}}
\] shown in p.252 of \cite{WW}.  (The convention of Bernoulli numbers in
\cite{WW} is different from us.)
 So $R_{\mathrm{osc}}$ satisfies the classical limit condition.

Also from the symmetry of the differential equation
we see that the orthogonal condition $R_{\mathrm{osc}}(-\hbar) ^t R_{\mathrm{osc}}(\hbar) =1$
since $R_{\mathrm{osc}}(-\hbar) ^t R_{\mathrm{osc}}(\hbar)=1$ mod $q$.

Since $\phi _{\ka,q}$ has the nonequivariant limit obviously,
$R_{\mathrm{osc}}$ has the nonequivariant limit of $\kl =0$.
Finally, we prove the equivariant homogeneity of
$R_{\mathrm{osc}}$. First notice that $f_q$ is quasi-homogeneous
of degree 1, if we assign degree 1 to the integration variables
and assign degree 2 (resp. degree 1) to $q_i$ (resp.
$\lambda_i$). Expand $f^\sigma _{q}$ at critical point $u_\sigma$
and so that after a coordinate change $f^\sigma _{q} = u_\sigma -
\sum y_j ^2$ for some local variables $y_j$. Now we see that
$$\int e^{f_q/\hbar} \omega = e^{u_i/\hbar}\int e^{-\sum _j
y_j^2/\hbar} (k_0 + k_1 y^2 + k_2 y^4...) \prod _j dy_j$$ where
$y_j$ has degree $1/2$ and $k_m$ has degree $-m-\frac{1}{2}\dim
X$.
Now using the asymptotic formula 
of the last integral and also of $\int e^{f_q/\hbar} \phi _{\ka ,
q}\omega$, we conclude that the degree of $R_k$ is $-k$. $\Box$

%

\bigskip

\subsection{The Virasoro constraints for flag manifolds}

In this section we explain in short
how to prove the Virasoro conjecture for flag
manifolds, applying Givental's theory \cite{Gh}.

Let $X$ be a complex projective algebraic manifold
with $H^{\mathrm{odd}}(X)=0$.
To study the so-called gravitational descendent Gromov-Witten invariants,
consider a generating function called
the genus $g$ descendent potential
$$F_X^g({\bf t}:=(t_0,t_1,...)) = \sum _{m, d\in H_2(X,\ZZ )}
\frac{Q^d}{m!} \int _{[\overline{M}_{g,m}(X,d)]^{virt}}
\prod _{i=1}^m(\sum_{k=0}^{\infty} \psi _i ^k ev_i ^* t_k  )$$
and the total descendent potential
$$Z_X = \exp \left(\sum _{g=0,1,..}\epsilon ^{g-1}F_X ^g \right).$$     Here
$t_i$ are in $H^*(X)$. The potential will be
considered as a formal function on
${\bf t}(\hbar ) = t_0 + t_1\hbar + t_2 \hbar ^2 + ...$.
Define $q$-coordinates  by the dilaton shift
${\bf q}(\hbar ) =q_0+q_1\hbar + q_2\hbar+...:= {\bf t}(\hbar ) -1\hbar$. Here
$1$ is the identity class in $H^*(X)$.

When $X$ is a point, Kontsevich proved
the Witten conjecture that the total descendent potential
$Z _{pt}$ is annihilated by specific
quadratic differential operators $\hat{L}_m + \delta _{m,0}/16$,
$m=-1,0,1,...$ with commuting
relation $[\hat{L}_m + \delta _{m,0}/16, \hat{L}_{m'} + \delta _{m',0}/16]
= (m-m') (\hat{L}_{m+m'} + \delta _{{m+m'},0}/16)$.
The commutation relation is the Lie algebra of
vector fields $-x^{m+1} \frac{d}{dx}$ on the line.
The first four of them are as follows.
\begin{eqnarray*}
 \hat{L}_{-1} & = & q_0^2 /2\epsilon + \sum _{m\ge 0} q_{m+1} \partial _m,\\
 \hat{L}_{0} &=& \sum _{m\ge 0} (m+1/2) q_m \partial _m, \\
 \hat{L}_1 &=& \epsilon \partial _0 ^2 /8 + \sum _{m\ge 0} (m+1/2)(m+3/2) q_m \partial _{m+1}, \\
 \hat{L}_2 &=& 3\epsilon \partial _0 \partial _1 /4 + \sum _{m\ge 0} (m+1/2) (m+3/2) (m+5/2) q_m \partial
 _{m+2}.
\end{eqnarray*}

Eguchi - Hori - Jinzenji - Xiong
and Katz \cite{E1, E2} extended the Witten conjecture
for Grassmannians $X$ and for all target spaces $X$, respectively. The extended conjecture is called the
Virasoro conjecture:
$(\mathcal{L}_m^X + \delta _{m,0}/16 ) Z_X =0$ where $\mathcal{L}_m^X$ are defined by data of
cohomology of $X$ and Chern classes of $X$, which will be  specified later in theorem \ref{Voperator}.
However $\mathcal{L}^X_{-1} =\sum _{\ka ,\kb }\frac{1}{2\epsilon} q_0^{\ka}q_0^{\kb} \eta _{\ka\kb}
+\sum_{m\ge 1,\ka }q_m^{\ka} \frac{\partial}{\partial q_{m-1}^a}$
and $\mathcal{L}^X_{-1}Z_X=0$ means the string equation.
Here $\eta _{\ka\kb}$ is the Poincare metric and ${\bf t}(\hbar)
=\sum _{\ka} t^{\ka}_0\phi _{\ka} + \sum _\ka t^{\ka}_1\phi _{\ka} \hbar+ ... $
with a fixed basis $\{\phi _{\ka} \}$ of $H^*(X)$.

If $H$ denotes the vector space $H^*(X)$ with Poincar\'e pairing
$(,)$, then the quotient ring $H((\hbar))$, of formal power
series of $\hbar$ over $H$, is endowed with a symplectic form
$\Omega$ defined by
$$\Omega (f,g)=\frac{1}{2\pi i} \oint (f(-\hbar),g(\hbar))d\hbar$$
for $f$ and $g$ in $H((\hbar ))$. So, $\mathcal{H}=H((\hbar))$ is
an infinite dimensional symplectic vector space. Given a
transformation $T$ of $\mathcal{H}$ which is infinitesimally
symplectic, that is $\Omega (Tf,g)+\Omega (f,Tg)=0$, define a
differential operator $\hat{T}$ as a quantization of $T$ as
follows. First, consider a quadratic function $\tilde{T}$
associated with $T$ by assignment $f\mapsto \frac{1}{2}\Omega
(f,Tf)$. Then take a quantization $\hat{T}$ of $\tilde{T}$ by the
rule: $p_i p_j \mapsto \epsilon \frac{\partial}{\partial
q_i}\frac{\partial}{\partial q_j}$, $p_iq_j \mapsto q_j
\frac{\partial}{\partial q_i}$, $q_iq_j \mapsto q_iq_j/\epsilon$
with Daboux coordinates $(p,q)$ of polarization
$\mathcal{H}=\mathcal{H}_++\mathcal{H}_-$ where $\mathcal{H}_+$
is the subspace of nonnegative power series of $\hbar$.  For
example, if $H=\CC ^2 $ with the standard inner product and $T$
is the multiplication operator by $1/\hbar$, i.e., $Tf=f/\hbar$
with $$f=...-(p_2^1 + p_2 ^2) /\hbar ^3 + (p_1^1 + p_1^2) /\hbar
^2 -(p_0^1 + p_0^2) /\hbar + (q_0^1 + q_0^2) + (q_1^1+q_2^2)
\hbar + ...$$ (here superscripts are indices), then
$$\hat{T} = \sum _{i=1,2} \left(\frac{(q_0 ^i)^2}{2\epsilon} +\sum _{m\ge 0} q_{m+1}^i \frac{\partial}{\partial q_m ^i}\right).$$
In fact, when $H=\CC$, and $D=\hbar\frac{d}{d\hbar }\hbar$, the
quantization of $D_m=\hbar ^{-1/2}D^{m+1}\hbar ^{-1/2} $ are
exactly the Virasoro operators $\hat{L}_m$ for $X=$point.

Define a transformation $S_t$ on $H[[\hbar ^{-1}]]$
by $(a,S_tb)= <<a,\frac{b}{\hbar -\psi}>>$, where
\begin{eqnarray*} &&
<<a,\frac{b}{\hbar - \psi }>>
= (a,b) + \\ &&
\sum _{0\ne d\in H_2(X),m=0,l=0}^{\infty, \infty}
\frac{Q^d}{ \hbar ^{l+1} m!}\int _{[\overline{M}_{0,m+2}(X,d)]^{virt}}
(ev_1^*a) (\prod _{i=2}^m ev^*_i t ) (ev^*_{m+2} b) \psi ^l _{m+2} .
\end{eqnarray*}

As introduced for $X=Fl(n+1)$, for general $X$ we have notions of
semi-simple quantum cohomology, canonical coordinates $u_{\ka}$,
a pencil of flat connections and an asymptotic fundamental
solution $\Psi R e^{U/\hbar}$. Here $R$ is form of $1+R_1\hbar
+...$ satisfying the orthogonality condition $R^*(-\hbar)R(\hbar
)=1$. Such $R$ is unique up to right multiplication by diagonal
matrices. However $R$ is uniquely determined if the homogeneity
condition is imposed. Now if $T$ is $S$ or $R(\hbar)$, then $T$
could be viewed as a symplectic transformation of suitable
completions of $\mathcal{H}$ since $T^*(\hbar)T(-\hbar)=1$. Let
$\hat{T}=\exp (\ln T)^{\hat{}}$. Denote $(q^1 (\hbar ),...,q^N
(\hbar ))=\Psi ^{-1} {\bf q}(\hbar ) $
 for ${\bf q}(\hbar ) \in \mathcal{H}_+$
and define  an operator $\hat{\Psi}$ by  $f(\Psi ^{-1}{\bf
q})\mapsto f({\bf q})$. Notice that $q^i (\hbar)$-coordinate
system is based on the orthonormal frame $\Delta _\ka
\frac{\partial}{\partial u_\ka}$.
 The homogeneity condition of
$R$ is $E (\Psi R e^{U/\hbar}) =\mu (\Psi R e^{U/\hbar })$, where
$E=\hbar \partial _{\hbar} + \sum u_{\ka} \partial _{u_{\ka}}$
and $\mu = \mathrm{diag} (\deg (\phi_1)-\dim_\CC X/2 , ...., \deg
(\phi_N) - \dim_\CC X /2 )$.
On the other hand $S$ satisfies the homogeneity condition of $E S
= \mu S + S (\mu + \rho /\hbar)^*$, where $\rho$ is the operator
of multiplication by $c_1(T_X)$ in ordinary cohomology ring.
Define
$$\mathcal{L}^{X}_m = \hat{S}^{-1}_u\hat{\Psi}\hat{R} \hat{L}_m
\hat{R}^{-1}\hat{\Psi}^{-1} \hat{S}_u$$  for $m=-1,0,1,2,...$.
where $L_m$ is taken as $D_m$ with $H=H^*(X)$.

The following theorem in \cite{Gh} explicitly shows that the
operator $\mathcal{L}^{X}_m$ is completely determined by
topological terms.
\begin{Thm}\label{Voperator} The operator
 $\mathcal{L}^X_m$ is $\hat{L}_{m}^{\mu, \rho}+ \frac{\delta _{m,0}}{4}\mathrm{tr}
 (\mu\mu^*)$, where $L_{m}^{\mu, \rho} =
 \hbar^\mu \hbar^{-\rho} L_m \hbar^\rho \hbar^{-\mu} = \hbar^{-1/2}( \hbar\frac{d}{d\hbar}-\mu
 \hbar + \rho )^{m+1} \hbar^{-1/2}, m\ge -1.$
\end{Thm}

The Virasoro operators  $\hat{L}_{m}^{\mu, \rho}+ \frac{\delta
_{m,0}}{4}\mathrm{tr}
 (\mu\mu^*)$ agree with the operators in \cite{DZ}.
 The theorem holds for a conformal semi-simple Frobenius manifolds.

Let a projective manifold $X$ have a Hamiltonian torus $T$ action
with isolated fixed points and let $u$ be a semisimple point of
$H^*_T(X)$. The previous potentials have the obvious equivariant
counterpart. The following is shown in \cite{Gsemi} and is
reformulated in \cite{Gh} as stated here.

\begin{Thm}\label{potential}
In the equivariant setting of Gromov-Witten theory,
 $$ Z_X^T ({\bf t}(\hbar ))
= e^{C(u)} \hat{S}_u ^{-1} \hat{\Psi}\hat{R} e^{(U /\hbar)^{\hat{} }}
\prod _{i=1}^N Z_{pt} (q^i (\hbar ))$$
if $R(z)$ is normalized by the classical limit condition in theorem \ref{Rconjforflag}
and $$C(u)=\frac12 \int ^u \sum _i( R_1 )_{ii} du^i$$
is defined up to addition of constant.
\end{Thm}

{\em Remark.} According to \cite{Gh} the right side of
the equation of the above theorem does not depend on the choices
of a semi-simple point $u$, even though each term may depend on
the choices.

\bigskip

The theorem shows $$(\mathcal{L}^{X,T}_m + N\delta _{m,0}/16
)Z_X^T=0,$$ where $\mathcal{L}^{X,T}_m$ is the equivariant
counterpart of $\mathcal{L}_m^X$.
\begin{Thm}\label{formula}
The total descendent potential $Z_X$ of  a flag manifold $X$
coincides with $e^{C(u)} \hat{S}_u ^{-1} \hat{\Psi}\hat{R} e^{(U
/\hbar)^{\hat{} }} \prod _{i=1}^N Z_{pt} (q^i (\hbar ))$.
\end{Thm}

{\em Proof.} Take the nonequivariant limit of the equation of
theorem \ref{potential} at a semi-simple point $u$ in $H^2_T(X)$.
The left side of the equation is specialized to the ordinary total
descendent potential. The limit of $R$ on the right side exits.
The limit is  the homogeneous ordinary $R$ due to  theorem
\ref{Rconjforflag} combined with  the uniqueness of homogeneous
$R|_{H^2_T(X)}$.
$\Box$

\bigskip

Combined with theorem \ref{Voperator}, the above theorem shows
the following corollary.
\begin{Cor}\label{virasoro}
The Virasoro conjecture for flag manifolds $X$ holds:
  $(\mathcal{L}_m^X + N\delta _{m,0}) Z_X =0$, $m\ge -1$.
\end{Cor}

\bigskip

{\em Acknowledgment.} B.K. would like to thank A. Givental,
D. van Straten for useful discussions and J.-H. Yang for informing
the existence of the paper \cite{St}. The authors also thank
J. Byeon for numerous discussions on oscillatory integrals.
B.K. thanks staffs in ESI for their warm hospitality while his visit
to the institute, being writing the paper.
D.J. is supported by KOSEF 2000-2-10100-002-3.
B.K. is supported by
KOSEF 1999-2-102-003-5 and R03-2001-00001.

\vspace{+10 pt}
\noindent
Department of Mathematics \\
Pohang University of Science and Technology \\
Pohang, 790-784\\
Republic of Korea \\
{\tt joe@euclid.postech.ac.kr \\
bumsig@postech.edu}

\end{document}